\documentclass[12pt]{amsart}
\usepackage{graphicx}
\usepackage{graphics}
\usepackage{amsmath}
\usepackage{amscd}
\usepackage{latexsym}
\usepackage{subfigure}
\usepackage{caption}
\usepackage{hyperref}

\begin{document}

\textwidth 5.9in
\textheight 7.9in

\evensidemargin .75in
\oddsidemargin .75in

\newtheorem{Thm}{Theorem}
\newtheorem{Lem}[Thm]{Lemma}
\newtheorem{Cor}[Thm]{Corollary}
\newtheorem{Prop}[Thm]{Proposition}
\newtheorem{Rm}{Remark}

\def\a{{\mathbb a}}
\def\C{{\mathbb C}}
\def\A{{\mathbb A}}
\def\B{{\mathbb B}}
\def\D{{\mathbb D}}
\def\E{{\mathbb E}}
\def\R{{\mathbb R}}
\def\P{{\mathbb P}}
\def\S{{\mathbb S}}
\def\Z{{\mathbb Z}}
\def\O{{\mathbb O}}
\def\H{{\mathbb H}}
\def\V{{\mathbb V}}
\def\Q{{\mathbb Q}}
\def\Cn{${\mathcal C}_n$}
\def\CM{\mathcal M}
\def\CG{\mathcal G}
\def\CH{\mathcal H}
\def\CT{\mathcal T}
\def\CF{\mathcal F}
\def\CA{\mathcal A}
\def\CB{\mathcal B}
\def\CD{\mathcal D}
\def\CP{\mathcal P}
\def\CS{\mathcal S}
\def\CZ{\mathcal Z}
\def\CE{\mathcal E}
\def\CL{\mathcal L}
\def\CV{\mathcal V}
\def\CW{\mathcal W}
\def\IC{\mathbb C}
\def\IF{\mathbb F}
\def\IK{\mathcal K}
\def\IL{\mathcal L}
\def\IP{\bf P}
\def\IR{\mathbb R}
\def\IZ{\mathbb Z}

\title{A  simple class of infinitely many absolutely exotic manifolds}
\author{Selman Akbulut}
\thanks{Partially supported by NSF grants DMS 0905917}
\address{Department  of Mathematics, Michigan State University,  MI, 48824}
\email{akbulut@math.msu.edu }
\date{\today}
\begin{abstract} 
We show that the smooth $4$-manifold $M$ obtained by attaching a $2$-handle to $B^4$ along a certain knot $K\subset \partial B^4$  admits infinitely many absolutely exotic copies $M_n$, $n=0,1,2..$, such that each copy $M_n$  is obtained by attaching $2$-handle to a fixed compact smooth contractible manifold $W$ along the iterates $f^{n}(c)$ of a knot $c\subset \partial W$ by a diffeomorphism  
$f:\partial W \to \partial W$. This generalizes the example in author's 1991 paper, which corresponds to $n=1$ case. 
\end{abstract}

\date{}
\maketitle

\setcounter{section}{-1}

\vspace{-.4in}

\section{Construction}

A relative exotic structure on a compact smooth $4$-manifold $M$ with boundary, is a self diffeomorphism  $f:\partial M\to \partial M$, which extends to a self homeomorphism of $M$, but does not extend to a self diffeomorphism of $M$. If $F:M\to M$ a homeomorphism extending $f$, then the pull-back smooth structure $M_{F}$ provides a relative exotic copy of $M$. We say that $N$ is an absolutely exotic copy of $M$, if it is homeomorphic but not diffeomorphic to $M$ (no condition on the boundary). The technique introduced in \cite{ar} turns relative exotic structures to absolute exotic structures. This is done by choosing an invertable cobordism $H$ with $\partial H=H_{-}\sqcup H_{+}$ and $H_{-}\approx \partial M$, and then gluing $H$ it to the boundary of $M$ in two different ways.  Then the manifolds $M'=M\cup_{Id}H$ and $M''=M\cup_{f}H$ become absolutely exotic copies of each other.  Applying this construction to a cork  $W$  produces an absolutely exotic copy of $W$; and when $W$ is an infinite order loose-cork (\cite{a2}) then we get infinitely many absolutely exotic copies of $W$. This construction results a boundary  $H_{+}$, which consists of hyperbolic manifolds glued along tori.  

\vspace{.1in}

Ideally we want to produce small $4$-manifolds with simple boundaries, admitting absolutely exotic copies (with the hope of capping boundaries to get small closed exotic manifolds). One such example is the cusp $C^4$ (\cite{a3}) with a Seifert fibered space boundary, which is obtained by attaching a $2$-handle to $B^4$ along the trefoil knot with $0$-framing. Performing knot surgeries $C\leadsto C_{K}$ along  the torus inside (by using different knots $K$) provides infinitely many absolutely exotic copies of $C$. An interesting open problem here is to find corks inside of $C$, such that twisting $C$ along them induce the exotic copies $C_{K}$ of $C$. Here we produce another example $M^4$, similar to $C^4$, which also has a Seifert fibered space boundary, and is obtained  by attaching a $2$-handle to $B^4$ along a slice knot with $-2$ framing. But from this we can construct infinitely many different exotic copies of M, each obtained by cork-twisting along an infinite order loose-cork $W\subset M$ (\cite{a2}), rather than a knot surgery to $M$, as in the case of $C$ above.

\vspace{.1in}

\begin{figure}[ht]  \begin{center} 
 \includegraphics[width=.19\textwidth]{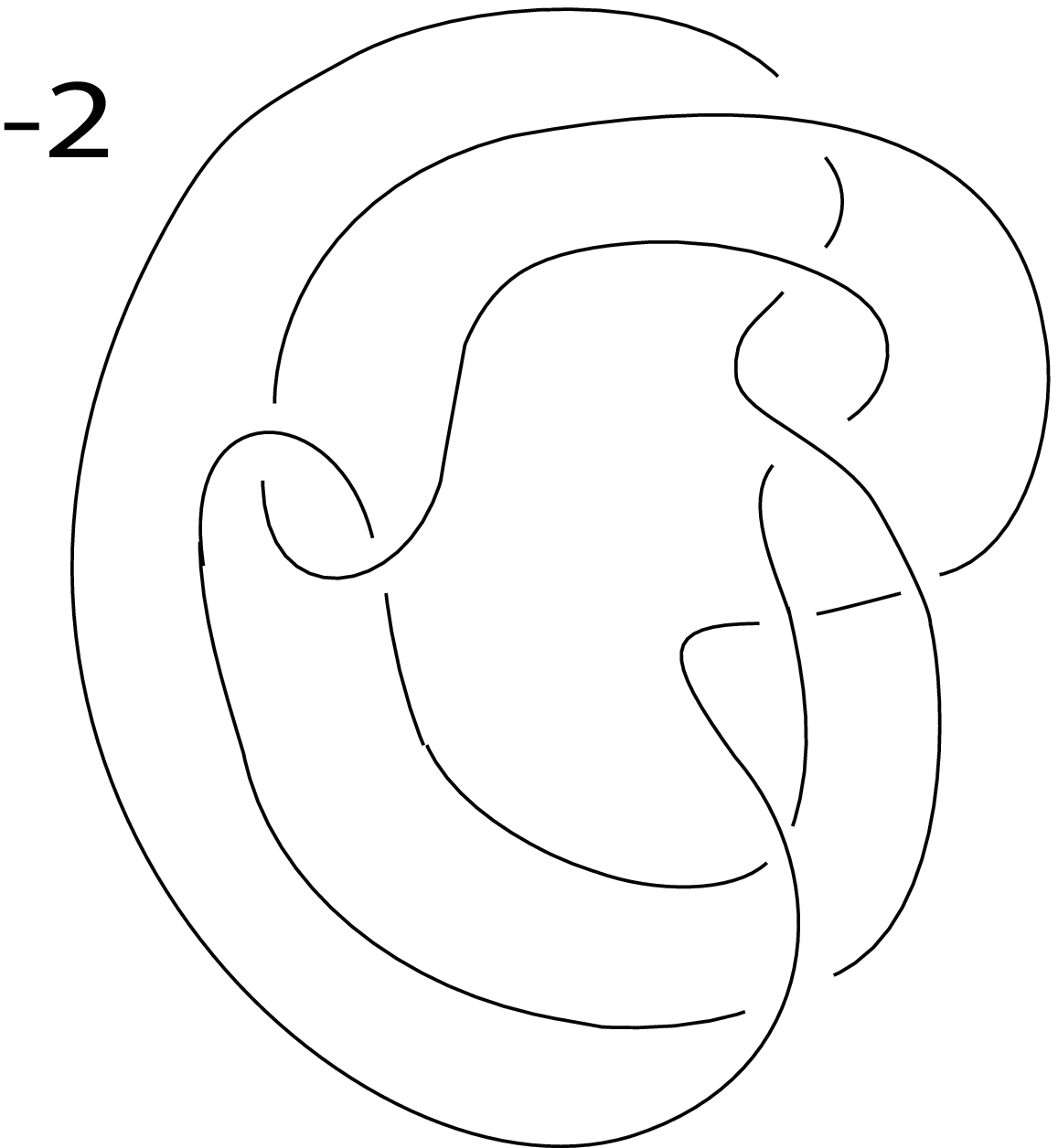}       
\caption{M}      \label{a1} 
\end{center}
 \end{figure}

\begin{Thm}\label{thm1}
The manifold $M$ of Figure~\ref{a1}, which is obtained by attaching $2$-handle to $B^4$ along the knot $K$ of Figure \ref{a1} with  $-2$ framing, admits infinitely many distinct absolutely exotic copies, and they can be detected twisting an infinite order loose-cork inside of $M$.
\end{Thm}

\begin{figure}[ht]  \begin{center}
 \includegraphics[width=.5\textwidth]{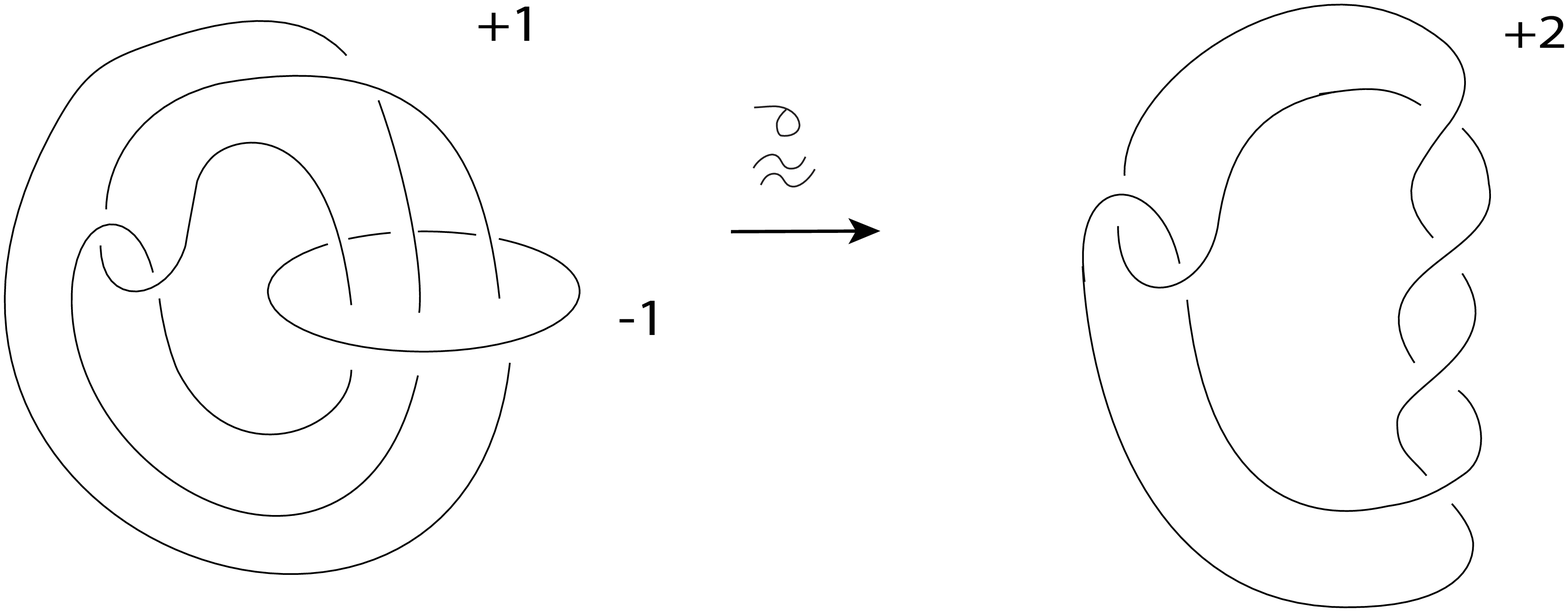}       
\caption{M}      \label{a2} 
\end{center}
 \end{figure}

\proof  First of all by blowing up and down as in Figure~\ref{a2}, we see that $\partial M$ can also be identified by $+2$ surgery to $(-4, 2)$ twist knot (stevedores knot). $\partial M$ is the small Seifert fibered space  $M (-2;1/2,3/4,7/9)$ (e.g. \cite{bw}, \cite{s}, \cite{t}), therefore its mapping class group is finite order.

\vspace{.1in}

\begin{figure}[ht]  \begin{center}
 \includegraphics[width=.7\textwidth]{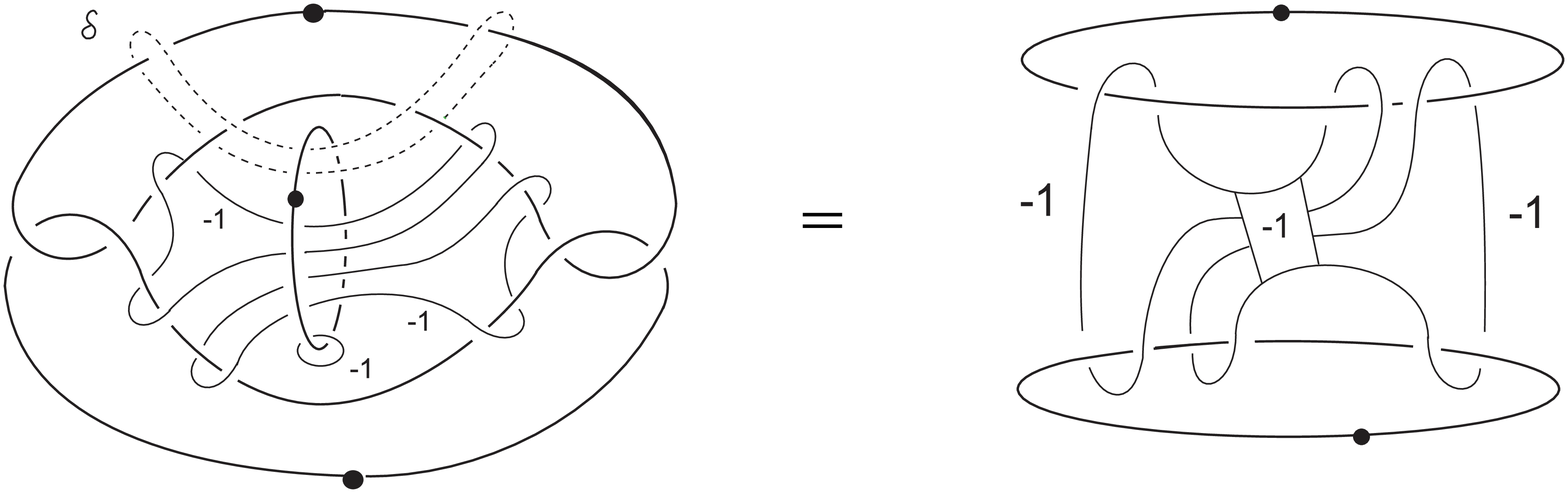}       
\caption{W}      \label{a3} 
\end{center}
 \end{figure}

Now, recall the infinite order loose-cork $W$ of \cite{a2}, which is shown in Figure~\ref{a3} (where two of its alternative handlebody pictures are given). The  order $n$ diffeomorphism $f_{n}: \partial W \to \partial W$ is obtained by a delta move along the curve $\delta \subset \partial W$ (the dotted curve in the figure). 

\begin{figure}[ht]  \begin{center}
 \includegraphics[width=.75\textwidth]{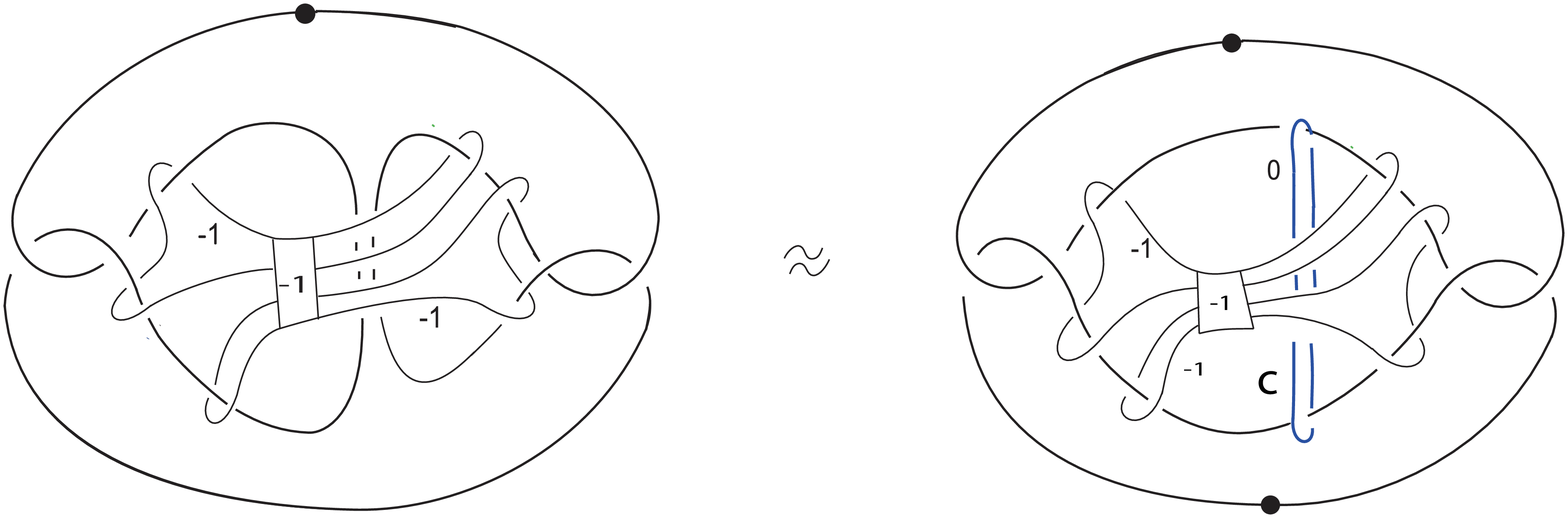}       
\caption{M}      \label{a4} 
\end{center}
 \end{figure}

Now check that the handlebody pictures of Figure~\ref{a4}  describe the manifold $M$ above (to see this cancel $1/2$- handle pairs). The second picture of Figure~\ref{a4} shows an imbedding $W\subset M$. That is, $M$ is obtained from $W$ by attaching a $2$-handle along the knot $c$ with $0$-framing. 

\begin{figure}[ht]  \begin{center}
 \includegraphics[width=.7\textwidth]{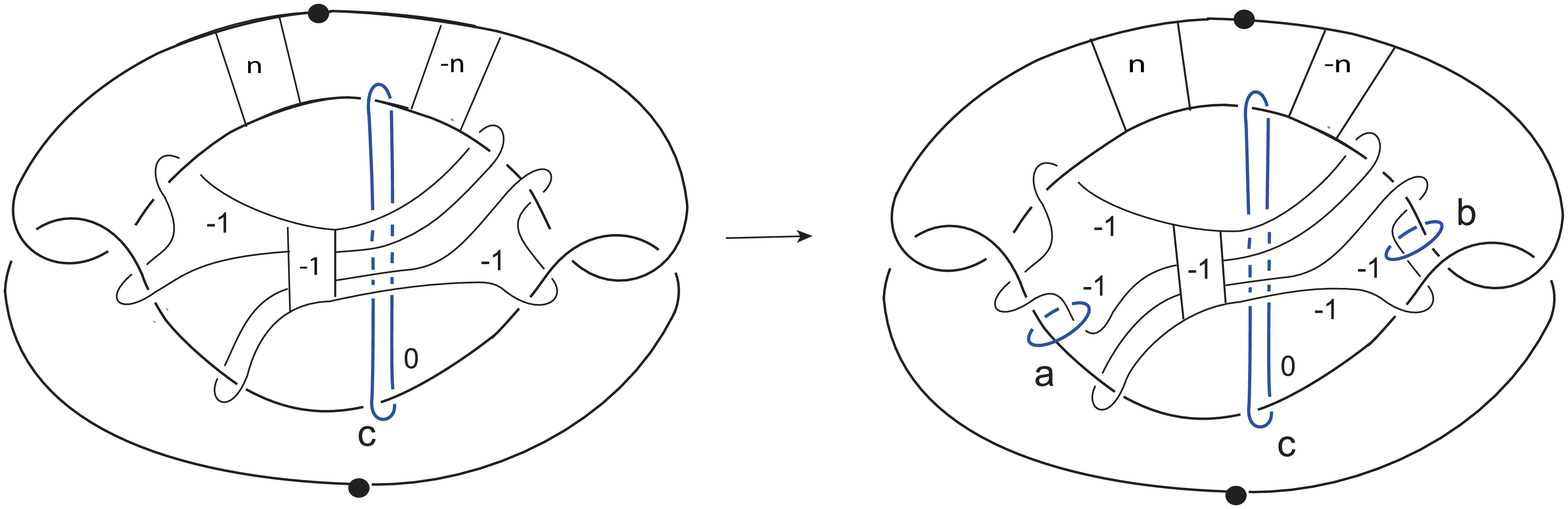}       
\caption{$M_{n} \leadsto S_{n}=M_{n}+a^{-1}+b^{-1}$}      \label{a5} 
\end{center}
 \end{figure}

Now apply $\delta$- moves to $W$ inside $M$,  $n$-times (where $\delta$ is chosen as in Figure~\ref{a3}), and call the resulting manifold $M_{n}$ (which is the first picture of Figure~\ref{a5}). We claim $\{M_{n}\}$ are exotic copies of $M$ rel boundary. To see this attach $2$-handles to $M_{n}$ along the knots $a$ and $b$ of the picture.

\begin{figure}[ht]  \begin{center}
 \includegraphics[width=.6\textwidth]{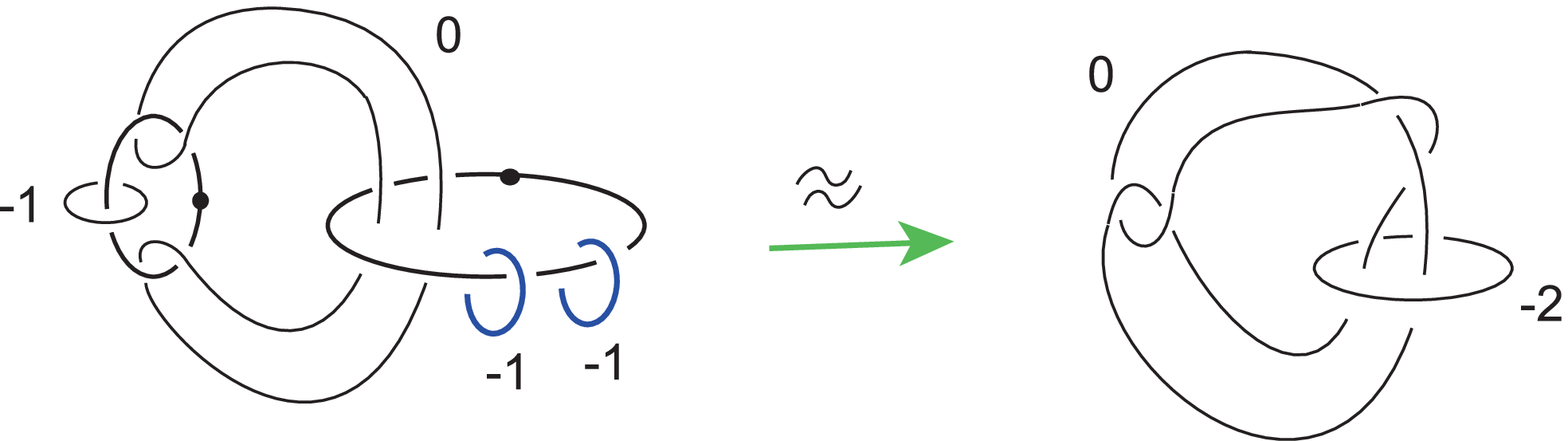}       
\caption{S}      \label{a7} 
\end{center}
 \end{figure}

Call the manifold  obtained from $M_{n}$ by attaching $2$-handles along $a$ and $b$ with $-1$ framings by $S_{n}= M_{n} +a^{-1}+b^{-1}$.  Now we proceed as in \cite{a2} by handle slides, to show that $S_{n}$ is the manifold obtained from the Stein manifold $S$ of Figure~\ref{a7} by the knot surgery using the twist knot $(-2, -n)$. Furthermore we can compactify $S$ into some closed symplectic manifold  $Z$ with $b_{2}^{+}(Z)>1$ (by \cite{lm},  \cite{ao}, or \cite{a4} p.108).  

\newpage

 \begin{figure}[ht]  \begin{center}
 \includegraphics[width=.7\textwidth]{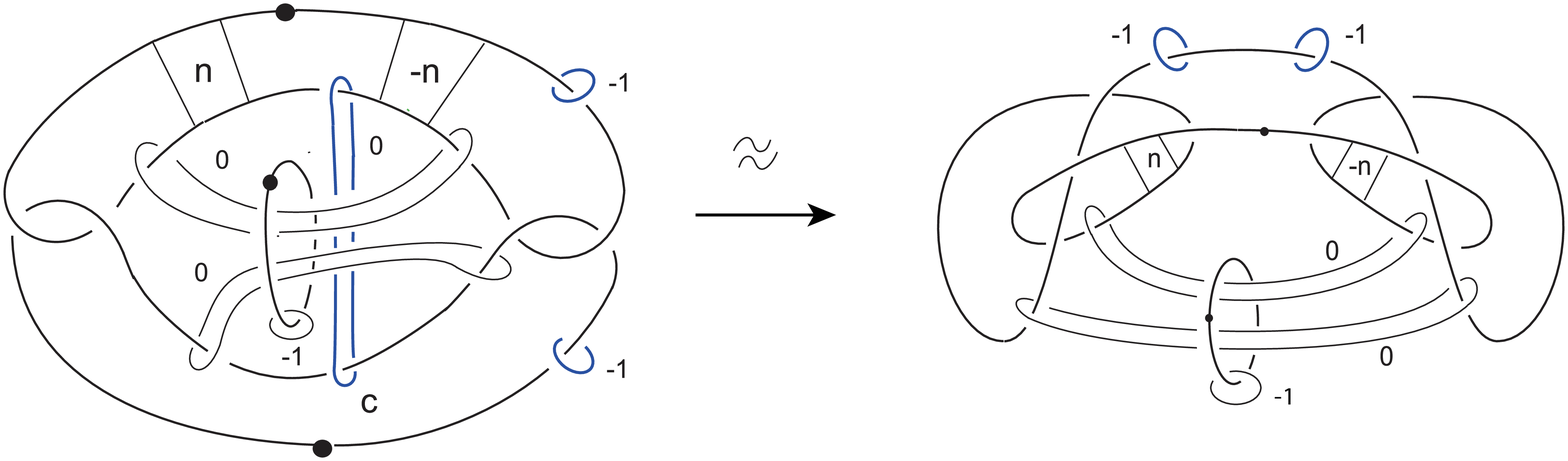}       
\caption{$S_{n}$}      \label{a6} 
\end{center}
 \end{figure}

 This shows that manifolds $\{M_{n}\}$  are exotic copies of $M$ rel boundary, and they are obtained by iterating $\delta$-moves to $f: \partial W \to \partial W $ inside $W\subset M$. Since the mapping class group of the Seifert fibered space $\partial M \approx M_{n}$ is finite \cite{bo}, by going to a subsequence we can assume that all $\{M_{n}\}$ are absolutely exotic copies of each other.
 
 \vspace{.1in} 
 
\begin{figure}[ht]  \begin{center}
 \includegraphics[width=.7\textwidth]{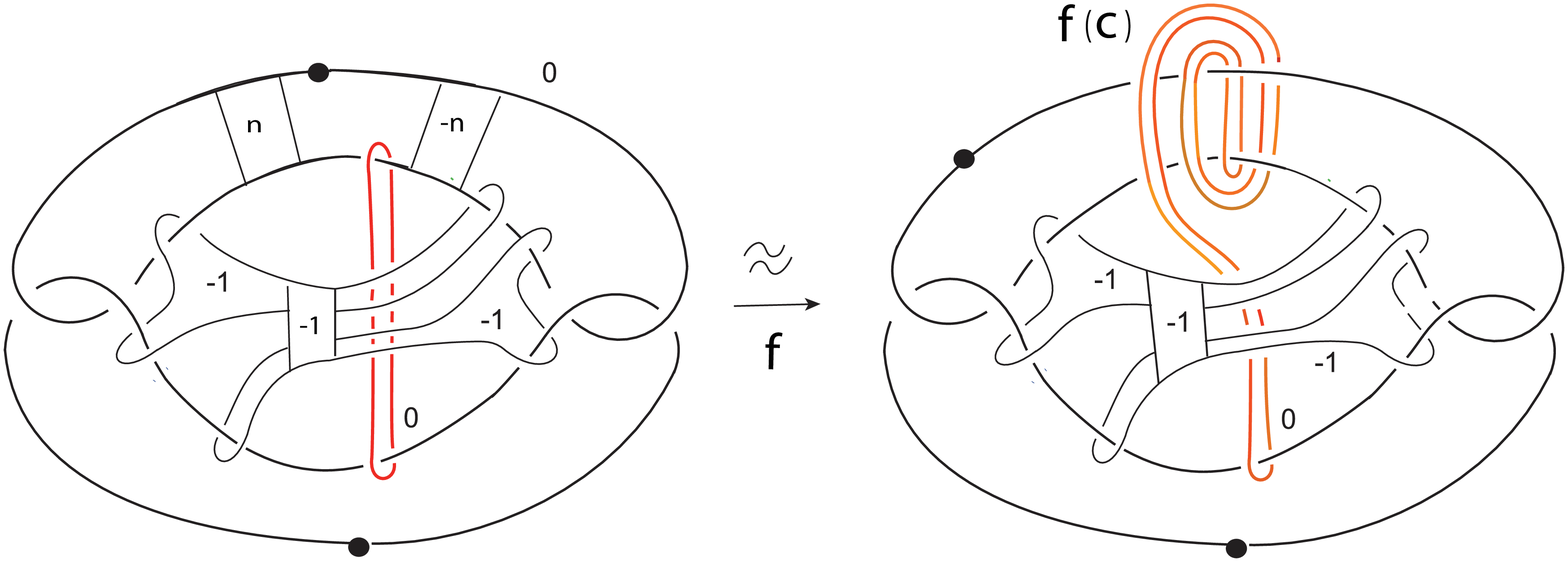}       
\caption{$M_n$}      \label{a8} 
\end{center}
 \end{figure}

 Now we analyze what an $n$-iterate of the $\delta$- move does to $M$, well it turns it into $M_{n}$, and a close inspection shows that $M_{n}$ is obtained from $M$ by attaching $2$-handle to $W$ along the loop $f^{n}(c)$ with $0$-framing, as shown in Figure~\ref{a8}.  Recall also, performing the $\delta$-move to $W$ inside of $M$, has the affect of attaching a cancelling pair of $2/3$-handles to $M$ and performing the diffeomorphism described in \cite{a1} resulting $M_{n}$.

 \begin{Rm}
 Note the new features of Theorem~\ref{thm1} which can not be reached by the techniques of \cite{ar}, they are:  (1) We don't need to modify the boundary of $M$ (by a homology cobordism) in order to construct its absolutely exotic copy.   (2) The construction here produces infinitely many absolutely exotic copies of $M$. (3) $M$ contains the tangent disc bundle of $S^2$ (an imbedded $-2$ sphere) and vice versa, so every smooth manifold which contains a $-2$ sphere contains a copy of  $M$ inside.  
\end{Rm}
 
 {Acknowledgements: I  thank Danny Ruberman for making suggestions and discussing some of the constructions of this paper with me.}


\begin{thebibliography}{99999}

 \bibitem[A1]{a1} S. Akbulut, {\em An exotic 4-manifold}, J. Differ. Geom. 33 (1991) 357--361.
 
 \bibitem[A2]{a2} S. Akbulut, {\em On infinite order corks}, PGGT 2016
 
\bibitem[A3]{a3}  S. Akbulut, {\em  A fake cusp and a fishtail.}   Turkish J. Math. 1 (1999) 19--31. 

  \bibitem[A4]{a4} S. Akbulut {\em $4$-Manifolds}, Oxford Graduate Texts in mathematics (2016) 
  ISBN: 9780198784869
  
   \bibitem[AO]{ao} S. Akbulut and B. Ozbagci {\em On the topology of compact Stein surfaces }, Int. Math. Res. Notices 15 (2002) 769-782
 

 \bibitem[AR]{ar} S. Akbulut and D.Ruberman {\em  An absolutely exotic contractible 4-manifold.} Commentarii Mathematici 
Helvetici, (2016) 91(1) 1-19


\bibitem[BO]{bo}  M. Boileau and J.-P. Otal {\em Scindements de Heegaard et groupe des homŽotopies des petites variŽtŽs de Seifert.}  Invent. Math.106 (1991), vol 1, 85 -107. 

\bibitem[BW]{bw}  Mark Brittenham and Ying-Qing Wu {\em 
The classification of exceptional Dehn surgeries on $2$-bridge knots}. Comm. Anal. Geom. 9 (2001), vol1, 97-113. 

\bibitem[S]{s} N. Saveliev {\em Invariants for homology $3$-spheres}, Encylopaedia of Math. Sci vol. 140, Springer-
Verlag ISBN 3-540-43796-7

\bibitem[LM]{lm} P. Lisca and G. Mattic {\em Tight contact structures and Seiberg-Witten invariants,} Invent. Math. 129 (1997) 509-525.

\bibitem[P]{p}  S. P. Plotnick {\em Vanishing of Whitehead groups of Seifert manifolds with infinite fundamental groups}, Comm. Math. Helv. 55 (1980) 654-667.

\bibitem[T]{t} B. Tosun {\em Tight small Seifert fibered manifolds with $e_{0}=-2$},\\ https://arxiv.org/pdf/1510.06948.pdf


 \end{thebibliography}
\end{document}